 \documentclass[final,1p,times]{elsarticle}
\usepackage{amssymb}
\usepackage{amsmath}
\usepackage{enumitem}
\usepackage{xcolor}
\begin{document}

\begin{frontmatter}

\title{Spectrally-accurate numerical method for acoustic scattering from doubly-periodic 3D multilayered media}
%\tnotetext[mytitlenote]{Fully documented templates are available in the elsarticle package on \href{http://www.ctan.org/tex-archive/macros/latex/contrib/elsarticle}{CTAN}.}

%% Group authors per affiliation:
\author{Min Hyung Cho\corref{mycorrespondingauthor}}
\address{Department of Mathematical Sciences, University of Massachusetts Lowell, Lowell, MA 01854}
\cortext[mycorrespondingauthor]{Corresponding author}
\ead{minhyung\_cho@uml.edu}

\begin{abstract}
A periodizing scheme and the method of fundamental solutions are used to solve acoustic wave scattering from doubly-periodic three-dimensional multilayered media. A scattered wave in a unit cell is represented by the sum of the near and distant contribution. The near contribution uses the free-space Green's function and its eight immediate neighbors. The contribution from the distant sources is expressed using proxy source points over a sphere surrounding the unit cell and its neighbors. The Rayleigh-Bloch radiation condition is applied to the top and bottom layers. Extra unknowns produced by the periodizing scheme in the linear system are eliminated using a Schur complement. The proposed numerical method avoids using singular quadratures and the quasi-periodic Green's function or complicated lattice sum techniques. Therefore, the proposed scheme is robust at all scattering parameters including Wood anomalies. The algorithm is also applicable to electromagnetic problems by using the dyadic Green's function.  Numerical examples with 10-digit accuracy are provided. Finally, reflection and transmission spectra are computed over a wide range of incident angles for device characterization. 
\end{abstract}


\end{frontmatter}

%\linenumbers

\section{Introduction}
Wave scattering from periodic structures and multilayered media plays a significant role in controlling the waves in modern electromagnetics and acoustic devices. Increasing numbers of  applications such as diffraction gratings, thin film photovoltaics \cite{atwater,kelzenberg}, photonic crystals \cite{joannopoulos2011photonic}, and meta-materials \cite{soukoulis2011past} utilize these structures to enhance the efficiency of devices. Thus, accurate and efficient numerical methods for their simulation are in very high demand. Traditionally, finite element methods (FEM) \cite{bao1995finite, monk2003finite, he2016spectral}, finite-difference time-domain (FDTD) methods \cite{taflovecomputational, haggblad2012consistent, chew19943d, winton2005fdtd}, and rigorous-coupled wave analysis (RCWA) or Fourier modal methods \cite{moharam1981rigorous, rokushima2006optics, li1996use, cho2008rigorous} are popular in physics and engineering fields due to their wide availability. However, most of the algorithms suffer from low accuracy and slow convergence in three dimensions. For example, low-order FEM suffers from the so-called pollution error or the accumulation of phase errors \cite{babuska1997pollution}. For a high frequency problem, degrees of freedom grow prohibitively large to maintain reasonable accuracy. FDTD has dispersion error and it can achieve only first or second order convergence. Both FDTD and FEM require some types of artificial boundary conditions at the truncation of the domain when dealing with the exterior domain problem. RCWA relies on intrinsically low-order staircase approximations of the layer interface. It is worth mentioning a recent development of a high-order perturbation of surface (HOPS) method \cite{nicholls2017numerical, hong2017high} that has shown promising results for shallow gratings in a low frequency regime.
 
 With the rapid improvement in computing power and fast algorithms such as the fast multipole method \cite{lapFMM, fmm1}, fast direct solver  \cite{mdirect, greengard2009fast}, and $\mathcal{H}$-matrix algorithms \cite{hackbusch}, Green's function or fundamental solution based methods have become more practical. Especially, for exterior domain problems, Green's function based methods have an advantage compared to other methods because Green's function naturally satisfies the outgoing radiation condition. Moreover, using the layer potentials or Green's second identity, the problem can be rewritten as a boundary or volume integral equation. The discretization of an integral equation using numerical quadrature rules results in a dense linear system that can be interpreted as the interaction between the source and target points. The dense linear system can then be solved by exploiting the low-rank structure of the matrix such as the fast multipole method  and fast direct solver. 
 \begin{figure}[t] %  figure placement: here, top, bottom, or page
   \centering
   \includegraphics[width=5in]{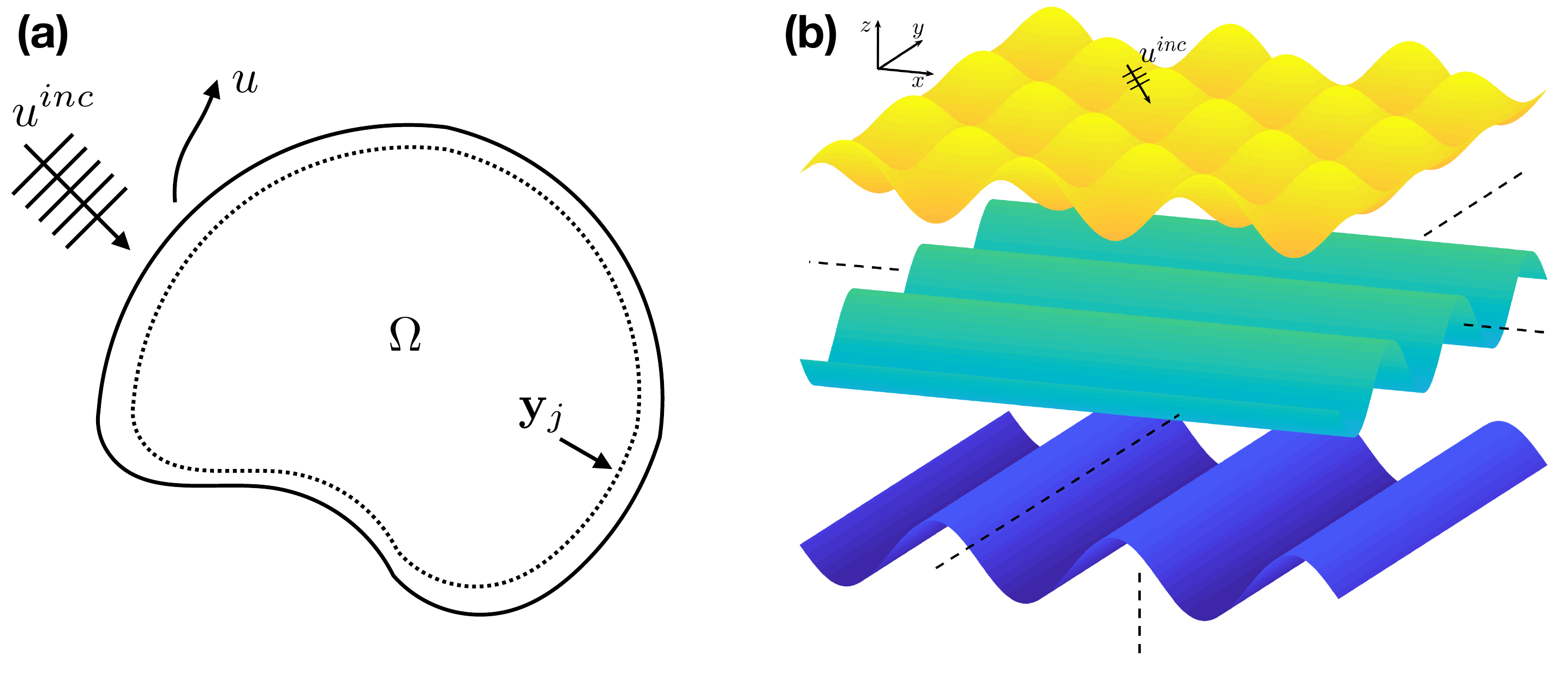}
   \caption{(a) Two-dimensional (2D) MFS for a sound-soft object $\Omega$ and (b) Problem set up for multilayered media with periodic shape (color online).}
   \label{fig:problem}
\end{figure}
 In this paper, an accurate and efficient numerical method based on a periodizing technique and the method of fundamental solutions (MFS) for multilayered media in three dimensions is proposed. The idea of MFS was first proposed by Kupradze and Aleksidze \cite{kupradze1964method}. The MFS is very similar to a boundary integral equation method where the source points have been moved away from the boundary. As a consequence, it avoids singular integrals at the cost of producing an ill-conditioned matrix. For example, consider a smooth sound-soft object $\Omega$ with the incident wave $u^{inc}$ in Fig. \ref{fig:problem}(a). The MFS represents the scattered field at $\mathbf{x} \in \mathbb{R}^2\backslash \Omega$ by a linear combination of the Green's function with unknown coefficients $c_j$
 \begin{align}
 u(\mathbf{x}) = \sum_{j=1}^N c_j G(\mathbf{x}, \mathbf{y}_j),\nonumber
 \end{align} where $G$ is the free-space Green's function of the Helmholtz equation and $\{\mathbf{y}_j\}_{j=1}^N$ is the artificial source points placed inside the domain. The Dirichlet boundary condition $u(\mathbf{x}_i)+u^{inc}(\mathbf{x}_i) = 0$ at the collocation point $\{\mathbf{x}_i\}_{i=1}^M$ on $\partial \Omega$ leads to an overdetermined linear system 
$$\sum_{j=1}^N c_j G(\mathbf{x}_i, \mathbf{y}_j) = -u^{inc}(\mathbf{x}_i)$$
and the best fit solution $c_j$ can be found. The interested readers are referred to a review article on MFS by Fairweather \cite{FaKa98}. In this paper, the problem consists of multiple smooth periodic surfaces. Figure \ref{fig:problem}(b) shows the problem setup. The plane wave that is quasi-periodic (periodic up to a phase) is incident in the top layer. The incident wave produces scattered waves that are also quasi-periodic. It is well known that a naive approach using the quasi-periodic Green's function faces many issues such as slow convergence and \textit{Wood anomalies} \cite{wood}. In order to overcome these issues, a new periodizing technique that uses only the free-space Green's function and contour integrals is introduced by Barnett and Greengard \cite{qplp,qpsc} for periodic objects. However, their extension to multilayered media was not straightforward. Therefore, boundary integral equation methods using the free-space Green’s function and auxiliary ring sources are applied to a large number of layers \cite{helmholtz_periodic} and many obstacles \cite{junlai} in two dimensions. A similar idea was applied to a three-dimensional (3D) Laplace equation \cite{gumerov}. However, for three dimensions, constructing an efficient and accurate quadrature rule for integral operators becomes challenging. Thus, MFS is used with the free-space Green's function for the unit cell and neighbors to represent near field contribution, and  proxy source points (or auxiliary sources) lying on a sphere that encloses the unit cell and its immediate neighboring cells are used to represent interaction from distant interfaces. Liu and Barnett \cite{liu2016efficient} applied a similar algorithm to doubly-periodic axisymmetric objects that are created by rotating a curve about the $z$-axis. In their approach,  due to axisymmetry, the three-dimensional problem reduces to a sequence of equations on the generating curves. The most common challenges and remedies for MFS methods are discussed very well in the reference. This work can be regarded as a generalization to 3D multilayered media that is not axisymmetric.   The proposed algorithm can be easily applied to Maxwell's equations using both electric and magnetic dyadic Green's functions. In Ref. \cite{kakulia2011method}, MFS and the spectral representation of quasi-periodic Green's function (which fails when $z=z'$) are applied to multilayered media.  Recently, the shifted Green's function with the domain decomposition method \cite{perez2018domain} is also applied for periodic multilayered media. Note that for planar-layered media, the layered media Green's function method is another effective approach because the layered media Green's function is constructed to satisfy interface conditions \cite{chew1999, cho2017efficient}. Consequently, most free-space methods can be used with minimal modification for an arbitrary distribution of objects embedded in layered media \cite{chen2016accurate, chen2018accurate}. However, computation of the layered media Green's function usually requires Sommerfeld integrals that need to be evaluated with high accuracy in a fast manner. There were many efforts to overcome these issues using window functions \cite{cai2000fast}, the wideband fast multipole method \cite{chocai12}, and the heterogeneous fast multipole method \cite{cho2017heterogeneous}. In Refs. \cite{chewbook, algorithmic_issue,caibook}, layered-media Green's function methods are reviewed. 

In the next two sections, a two-layer structure with the Dirichlet boundary condition is presented to illustrate the proposed method. Then, in Section 4, the method is extended to multilayered media with transmission boundary conditions. Numerical solutions are presented in Section 5. Finally, the paper concludes with a summary and future direction.

%% \linenumbers
%% main text
\begin{figure}[t] %  figure placement: here, top, bottom, or page
   \centering
   \includegraphics[width=5.2in]{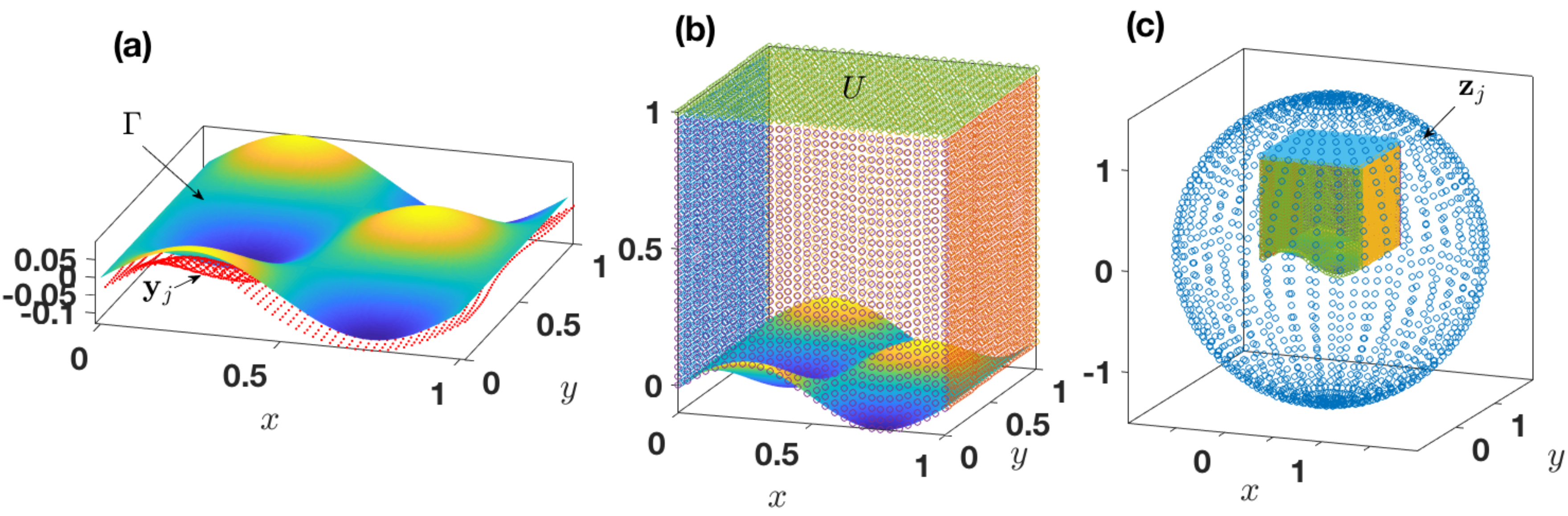} %low resolution
   \caption{(a) Layer interface ($\Gamma$) of the unit cell and MFS source points $\mathbf{y}_j$ (red dots) for the Dirichlet problem. (b) Left ($\Gamma_L$), right ($\Gamma_R$),  back ($\Gamma_B$), and front ($\Gamma_F$) walls, and artificial top layer ($U$). (c) The proxy source points $\mathbf{z}_j$ (color online).}
   \label{fig:dirichlet_domain}
\end{figure}
\section{MFS for a periodic interface with Dirichlet boundary condition}
In this section, a half-space is considered.  Let $\Gamma = \{\mathbf{x} = (x,y,z) | z = g(x,y), x \in [0, e_x], y\in [0, e_y]\}$ be the interface depicted in Fig. \ref{fig:dirichlet_domain}(a). The left, right, back, and front walls surrounding the unit cell are denoted by $\Gamma_L$, $\Gamma_R$, $\Gamma_B$, and $\Gamma_F$, respectively (See Fig. \ref{fig:dirichlet_domain}(b)). The Rayleigh-Bloch radiation condition will be enforced on the fictitious interfaces $U$ at $z = z_u$ in Fig. \ref{fig:dirichlet_domain}(b). 
%$\Gamma_L = \{\mathbf{x}_L = (x,y,z)| x = 0, y \in [0, e_y], z > g(0, y)\}$, $\Gamma_R = \{\mathbf{x}_R =(x,y,z)| x = e_x, y \in [0, e_y], z > g(e_x, y)\}$, $\Gamma_B = \{\mathbf{x}_B =(x,y,z)| x \in [0, e_x], y=0, z > g(x, 0)\}$, and $\Gamma_F = \{\mathbf{x}_F =(x,y,z)| x \in [0, e_x], y = e_y, z > g(x, e_y)\}$, respectively. 
Let $\Omega = \{(x,y,z)|  g(x,y) < z < z_u, x \in [0, e_x], y\in [0, e_y] \}$ be the region above $\Gamma$ and below $U$. The plane wave $u^{inc} = e^{i\mathbf{k}\cdot \mathbf{x}}$ is incident in $\Omega$ with the wavevector $\mathbf{k} = (k_x, k_y, k_z) = (k\sin{\phi^{inc}}\cos{\theta^{inc}}, k\sin{\phi^{inc}}\sin{\theta^{inc}}, k\cos{\phi^{inc}})$, where  $0 \leq  \theta^{inc}  < 2\pi$,  $\pi/2  < \phi^{inc} \leq \pi$, and $k = |\mathbf{k}|$. The incident wave is periodic up to a phase or quasi-periodic, namely,
\begin{align}
\alpha_x^{-1} u^{inc}(x+e_x, y, z) = \alpha_y^{-1} u^{inc}(x, y+e_y, z) =  u^{inc}(x,y,z),
\end{align}
where Bloch phases are defined by 
\begin{align}
\alpha_x = e^{ie_x k_x} \mbox{ and } \alpha_y = e^{ie_y k_y}.
\end{align} 
Then, the incident wave generates the scattered wave $u$. It is well known that the scattered wave $u$ satisfies the Helmholtz equation in the upper half-space, and it is quasi-periodic as well. Therefore, the boundary value problem for $u$ with the Dirichlet boundary condition can be written as
\begin{align}
&\Delta u+k^2 u = 0,\mathbf{x}\in\Omega,\\
&u+u^{inc} = 0, \mathbf{x}\in\Gamma,\\
&\left. u\right\vert_{\Gamma_L} =\alpha_x^{-1} \left. u\right\vert_{\Gamma_R}, \left. u\right\vert_{\Gamma_B} =\alpha_y^{-1} \left. u\right\vert_{\Gamma_F}, \\
&\left.  \frac{\partial u}{\partial \mathbf{n}}\right\vert_{\Gamma_L} = \alpha_x^{-1} \left. \frac{\partial u}{\partial \mathbf{n}}\right\rvert_{\Gamma_R}, \left.  \frac{\partial u}{\partial \mathbf{n}}\right\vert_{\Gamma_B} = \alpha_y^{-1} \left. \frac{\partial u}{\partial \mathbf{n}}\right\rvert_{\Gamma_F},
\end{align}
with the upward Rayleigh-Bloch radiation condition
\begin{align}
u(\mathbf{x}) = \sum_{m,n \in Z} a_{mn}^ue^{[i(\kappa_x^m x +\kappa_y^n y +k_u^{(m,n)}(z-z_u))]},~~ z \geq z_u, \label{RB_solution}
\end{align}
where $\kappa_x^m = k_x+2\pi m/e_x$, $\kappa_y^n = k_y+2\pi n/e_y$, and $k_u^{(m,n)} = \sqrt{k^2-(\kappa_x^m)^2 - (\kappa_y^n)^2}$ and the sign of the square root is taken as positive real or positive imaginary.  The coefficients $a^u_{mn}$ are the Bragg diffraction amplitudes of the reflected wave. The solution to this problem can be represented with
MFS using the quasi-periodic Green's function,
\begin{align}
u(\mathbf{x}) = \sum_{j=1}^N c_j G^{QP}_k(\mathbf{x}, \mathbf{y}_j), \label{MFS_original} 
\end{align}
where $\mathbf{y}_j$ is the artificial source points placed under $\Gamma$ toward the normal direction from the interface (red dots in Fig. \ref{fig:dirichlet_domain}(a)). The quasi-periodic Green's function $G^{QP}_k (\mathbf{x}, \mathbf{y})$ for the 3D Helmholtz equation is defined by
\begin{align}
G^{QP}_k(\mathbf{x}, \mathbf{y}) =  \sum_{m=-\infty}^{\infty} \sum_{n=-\infty}^{\infty}\alpha_x^m \alpha_y^n G_k(\mathbf{x}, \mathbf{y}+m\mathbf{e}_x+n\mathbf{e}_y),
\end{align}
where $\mathbf{e}_x = (e_x, 0, 0)$, $\mathbf{e}_y = (0, e_y, 0)$, and $G_k$ is the free-space Green's function of 3D Helmholtz equation with the wavenumber $k$ given by
\begin{align}
G_k(\mathbf{x}, \mathbf{y}) = \frac{e^{ik|\mathbf{x}-\mathbf{y}|}}{4\pi |\mathbf{x}-\mathbf{y}|}.
\end{align}
 %Here, the distance parameter for source points plays an crucial role determining the accuracy. Three methods are proposed in Ref. \cite{liu2016efficient}. 

%Then, the Dirichlet boundary condition for $\mathbf{x}\in\Gamma$ yields 
%\begin{align}
%u(\mathbf{x})  = \sum_{j=1}^N c_j G^{QP}_k(\mathbf{x}, \mathbf{y}_j) = -u^{inc}(\mathbf{x}), \mathbf{x} \in \Gamma. 
%\end{align}
%By discretizing the surface $\Gamma$ by $\{\mathbf{x}_i\}_{i=1}^{M}$, a rectangular linear system of equation for $\{c_j\}_{j=1}^{N}$
%\begin{align}
% \sum_{j=1}^N c_j G^{QP}_k(\mathbf{x}_i, \mathbf{y}_j) = -u^{inc}(\mathbf{x}_i),~ \mathbf{x}_i \in \Gamma, i=1,2,\cdots M 
%\end{align}
%can be obtained. 
However, in this approach, the quasi-periodic Green's function suffers from slow convergence and becomes impractical for 3D problems. There are many methods to remedy the slow convergence such as the Ewald summation method \cite{ewald1921berechnung, jordan1986efficient, arens2011analysing}, spatial-spectral splitting \cite{jorgenson1990efficient}, or a lattice sum \cite{lintonrev,otani2008periodic, denlinger2017fast}. Moreover, the quasi-periodic Green's function does not exist at Wood anomalies. The Wood anomalies are a set of special scattering parameters $\alpha_x$, $\alpha_y$, and $k$, for which the sum in the quasi-periodic Green's function diverges even if the problem is well-posed \cite{shipmanreview}. Non-robustness at these parameters is addressed by replacing the quasi-periodic Green's function with the quasi-periodic Green's function for other boundary conditions \cite{CWnystrom,horoshenkov,brunoqp3d,bruno2014} or using the free-space Green's function with immediate neighbors and equivalent sources representation for far field contributions. In the next sections, a periodizing scheme using the second approach is presented.

\section{Periodizing scheme for a periodic layer interface with the Dirichlet boundary condition}
In this section, a periodizing algorithm for the half-space (single interface) is presented using MFS with the finite sum of the free-space Green's function and proxy source points over a sphere. A similar algorithm with the spherical harmonic expansion (instead of proxy source points) is used for acoustic scattering from doubly-periodic axisymmetric obstacles in three dimensions \cite{liu2016efficient}. 

The scattered field $u$ is decomposed by near interaction with the unit cell and its surrounding eight immediate neighboring cells and far field contribution using proxy source points over a sphere as
\begin{align}
u(\mathbf{x}) = \sum_{j=1}^N  \sum_{m=-1}^1 \sum_{n=-1}^1c_j \alpha_x^m \alpha_y^n G_k(\mathbf{x}, \mathbf{y}_j+m\mathbf{e}_x+n\mathbf{e}_y)+\sum_{j=1}^Ld_jG_k(\mathbf{x},\mathbf{z}_j) , \label{MFS_periodic} 
\end{align}where $\mathbf{y}_j$ are the MFS source points placed below the layer interface (See Fig. \ref{fig:dirichlet_domain}(a)) and $\mathbf{z}_j$ are the proxy source points over a sphere enclosing the unit cell and its immediate neighbors (See Fig. \ref{fig:dirichlet_domain}(c)). These can be interpreted as equivalent source points that represent effects from distant copies of the unit cell. In the following, a linear system for the MFS coefficients $\mathbf{c}$, the proxy strength unknowns $\mathbf{d}$, and the Bragg coefficients $\mathbf{a}$ is obtained  by imposing (a) the Dirichlet boundary condition at $\Gamma$, (b) the quasi-periodicity on the surrounding walls $\Gamma_L$, $\Gamma_R$, $\Gamma_B$, and $\Gamma_F$, and (c) the upward Rayleigh-Bloch radiation condition at the fictitious layer $U$:

\begin{enumerate}[label=(\alph*)]
\item Dirichlet boundary condition ($u(\mathbf{x}) = -u^{inc}(\mathbf{x}) $, $\mathbf{x} \in \Gamma$)\\
The scattered field $u$ must satisfy the Dirichlet boundary condition. Thus, for uniformly distributed $\{\mathbf{x}_i\}_{i=1}^M \in \Gamma$,
\begin{align}
u(\mathbf{x}_i)=\sum_{j=1}^N  \sum_{m=-1}^1 \sum_{n=-1}^1c_j \alpha_x^m \alpha_y^n G_k(\mathbf{x}_i, \mathbf{y}_j+m\mathbf{e}_x+n\mathbf{e}_y)+\sum_{j=1}^Ld_jG_k(\mathbf{x}_i,\mathbf{z}_j) = -u^{inc}(\mathbf{x}_i)
\end{align}
or in a matrix form
\begin{align}
\mathbf{A}\mathbf{c}+\mathbf{B}\mathbf{d} = \mathbf{f}, \label{eq:dirichlet}
\end{align}
where
\begin{align}
\mathbf{A}& =\left[ \sum_{m=-1}^1 \sum_{n=-1}^1\alpha_x^m \alpha_y^n G_k(\mathbf{x}_i, \mathbf{y}_j+m\mathbf{e}_x+n\mathbf{e}_y) \right], \mathbf{f} = \left[-u^{inc}(\mathbf{x}_i) \right],  i=1,2,\cdots M, j=1,2,\cdots N, \nonumber\\
\mathbf{B} &=\left[G_k(\mathbf{x}_i, \mathbf{z}_j)\right],   i=1,2,\cdots M, j=1,2,\cdots L,\nonumber\\
\mathbf{c} &=[c_j], j=1,2,\cdots N, \mathbf{d} = [d_j], j=1,2,\cdots L.
\end{align}

\item Quasi-periodic boundary conditions\\
The scattered field must satisfy the quasi-periodic boundary conditions
\begin{align}
\left. u\right\vert_{\Gamma_L} - \alpha_x^{-1} \left. u\right\vert_{\Gamma_R} = 0,& \left.  \frac{\partial u}{\partial \mathbf{n}}\right\vert_{\Gamma_L} - \alpha_x^{-1} \left. \frac{\partial u}{\partial \mathbf{n}}\right\rvert_{\Gamma_R} = 0,\\
\left. u\right\vert_{\Gamma_B} - \alpha_y^{-1} \left. u\right\vert_{\Gamma_F} = 0,& \left. \frac{\partial u}{\partial \mathbf{n}}\right\vert_{\Gamma_B} - \alpha_y^{-1} \left. \frac{\partial u}{\partial \mathbf{n}} \right\rvert_{\Gamma_F} = 0.
\end{align}
For the sake of simplicity, only the quasi-periodic boundary condition on the left ($\Gamma_L$) and right ($\Gamma_R$) walls is presented.  The scattered fields at $\mathbf{x}_L \in \Gamma_L$ and $\mathbf{x}_R \in \Gamma_R$ are
\begin{align}
u(\mathbf{x}_L) = \sum_{j=1}^N  \sum_{m=-1}^1 \sum_{n=-1}^1c_j \alpha_x^m \alpha_y^n G_k(\mathbf{x}_L, \mathbf{y}_j+m\mathbf{e}_x+n\mathbf{e}_y)+\sum_{j=1}^Ld_jG_k(\mathbf{x}_L,\mathbf{z}_j), \\
 u(\mathbf{x}_R) =  \sum_{j=1}^N  \sum_{m=-1}^1 \sum_{n=-1}^1c_j \alpha_x^{m} \alpha_y^n G_k(\mathbf{x}_R, \mathbf{y}_j+m\mathbf{e}_x+n\mathbf{e}_y)+\sum_{j=1}^Ld_jG_k(\mathbf{x}_R,\mathbf{z}_j).
 \end{align}
Then, the quasi-periodic boundary condition $\left. u\right\vert_{\Gamma_L} - \alpha_x^{-1} \left. u\right\vert_{\Gamma_R} = 0$  and the translational symmetry \cite{qpsc} result in
\begin{align}
 \sum_{j=1}^N  c_j  \left( \alpha_x^{-2} \sum_{n=-1}^1\alpha_y^n G_k(\mathbf{x}_R, \mathbf{y}_j-\mathbf{e}_x+n\mathbf{e}_y) -\alpha_x^{1}  \sum_{n=-1}^1 \alpha_y^n G_k(\mathbf{x}_L, \mathbf{y}_j+\mathbf{e}_x+n\mathbf{e}_y) \right)\nonumber\\
+\sum_{j=1}^Ld_j \left(\alpha_x^{-1}G_k(\mathbf{x}_R,\mathbf{z}_j)-G_k(\mathbf{x}_L,\mathbf{z}_j) \right)= 0.
\end{align}
All other conditions yield very similar equations. Therefore, the quasi-periodic conditions produce a linear system for $\mathbf{c}$ and $\mathbf{d}$ as follows:
\begin{align}
\mathbf{P}\mathbf{c}+\mathbf{Q}\mathbf{d} = 0,\label{eq:qbc}
\end{align}
where
\begin{align}
\mathbf{P} = 
\left[\begin{array}{c}\alpha_x^{-2} \displaystyle \sum_{n=-1}^1\alpha_y^n G_k(\mathbf{x}_R, \mathbf{y}_j-\mathbf{e}_x+n\mathbf{e}_y) -\alpha_x^{1}  \sum_{n=-1}^1 \alpha_y^n G_k(\mathbf{x}_L, \mathbf{y}_j+\mathbf{e}_x+n\mathbf{e}_y) \\
\alpha_x^{-2} \displaystyle\sum_{n=-1}^1\alpha_y^n \frac{\partial G_k}{\partial \mathbf{n}}(\mathbf{x}_R, \mathbf{y}_j-\mathbf{e}_x+n\mathbf{e}_y) -\alpha_x^{1}  \sum_{n=-1}^1 \alpha_y^n \frac{\partial G_k}{\partial \mathbf{n}}(\mathbf{x}_L, \mathbf{y}_j+\mathbf{e}_x+n\mathbf{e}_y) \\
\alpha_y^{-2} \displaystyle\sum_{m=-1}^1\alpha_x^m G_k(\mathbf{x}_F, \mathbf{y}_j+m\mathbf{e}_x-\mathbf{e}_y) -\alpha_y^{1}  \sum_{m=-1}^1 \alpha_x^m G_k(\mathbf{x}_B, \mathbf{y}_j+m\mathbf{e}_x+\mathbf{e}_y) \\
\alpha_y^{-2}\displaystyle \sum_{m=-1}^1\alpha_x^m \frac{\partial G_k}{\partial \mathbf{n}}(\mathbf{x}_F, \mathbf{y}_j+m\mathbf{e}_x-\mathbf{e}_y) -\alpha_y^{1}  \sum_{m=-1}^1 \alpha_x^m \frac{\partial G_k}{\partial \mathbf{n}}(\mathbf{x}_B, \mathbf{y}_j+m\mathbf{e}_x+\mathbf{e}_y)\end{array}\right]
\end{align}
and
\begin{align}
\mathbf{Q} = \left[\begin{array}{c}
\displaystyle \alpha_x^{-1}G_k(\mathbf{x}_R,\mathbf{z}_j)-G_k(\mathbf{x}_L,\mathbf{z}_j)  \\
\displaystyle \alpha_x^{-1}\frac{\partial G_k}{\partial \mathbf{n}}(\mathbf{x}_R,\mathbf{z}_j)-\frac{\partial G_k}{\partial \mathbf{n}}(\mathbf{x}_L,\mathbf{z}_j)  \\
\displaystyle \alpha_y^{-1}G_k(\mathbf{x}_F,\mathbf{z}_j)-G_k(\mathbf{x}_B,\mathbf{z}_j)  \\
\displaystyle\alpha_y^{-1}\frac{\partial G_k}{\partial \mathbf{n}}(\mathbf{x}_F,\mathbf{z}_j)-\frac{\partial G_k}{\partial \mathbf{n}}(\mathbf{x}_B,\mathbf{z}_j) \end{array}\right].
\end{align}

\item Upward Rayleigh-Bloch radiation condition\\
The upward radiation condition is imposed at the artificial boundary $U$ using the Rayleigh-Bloch expansion. For uniformly distributed $\{ \mathbf{x}_i = (x_i, y_i, z_u) \}_{i=1}^{N_t} \in U$, the MFS representation and its normal derivative must be equal to the Rayleigh-Bloch expansion and its normal derivative, respectively. Therefore,
\begin{align}
& \sum_{j=1}^N  \sum_{m=-1}^1 \sum_{n=-1}^1c_j \alpha_x^m \alpha_y^n G_k(\mathbf{x}_i, \mathbf{y}_j+m\mathbf{e}_x+n\mathbf{e}_y)+\sum_{j=1}^Ld_jG_k(\mathbf{x}_i,\mathbf{z}_j) \nonumber\\
&~~~~~~~~~~~~~~~~~~~~~~~~~~~~~~~~~~~~~~~~~~~~~~~~~~~~~~~~~~~~~~-  \sum_{m=-R}^{R} \sum_{n = -R}^R a^u_{mn}e^{i(\kappa_x^m x_i+\kappa_y^n y_i)  } = 0,\\
 &\sum_{j=1}^N  \sum_{m=-1}^1 \sum_{n=-1}^1c_j \alpha_x^m \alpha_y^n \frac{\partial G_k}{\partial \mathbf{n}}(\mathbf{x}_i, \mathbf{y}_j+m\mathbf{e}_x+n\mathbf{e}_y)+\sum_{j=1}^Ld_j\frac{\partial G_k}{\partial \mathbf{n}}(\mathbf{x}_i,\mathbf{z}_j) \nonumber\\
&~~~~~~~~~~~~~~~~~~~~~~~~~~~~~~~~~~~~~~~~~~~~~~~~~~~~~~~~~~~~~~  +  \sum_{m=-R}^{R} \sum_{n = -R}^R a^u_{mn}ik_{u}^{(m,n)}e^{i(\kappa_x^m x_i+\kappa_y^n y_i)  } = 0,
 \end{align}
or in a matrix form
\begin{align}
Z\mathbf{c}+V\mathbf{d}-W\mathbf{a} = 0,\label{eq:radiation}
\end{align}
where
\begin{align}
Z &= \left[\begin{array}{c}
\displaystyle\sum_{m=-1}^1 \sum_{n=-1}^1\alpha_x^m \alpha_y^n  G_k(\mathbf{x}_i, \mathbf{y}_j+m\mathbf{e}_x+n\mathbf{e}_y) \\
\displaystyle\sum_{m=-1}^1 \sum_{n=-1}^1\alpha_x^m \alpha_y^n  \frac{\partial G_k}{\partial \mathbf{n}}(\mathbf{x}_i, \mathbf{y}_j+m\mathbf{e}_x+n\mathbf{e}_y)\end{array}\right], i=1,2,\cdots, N_t, j=1,2,\cdots N,\\
 V& = \left[\begin{array}{c}G_k(\mathbf{x}_i,\mathbf{z}_j) \\   \displaystyle  \frac{\partial G_k}{\partial \mathbf{n}}(\mathbf{x}_i,\mathbf{z}_j)\end{array}\right], i=1,2,\cdots, N_t, j=1,2,\cdots L,\\
W &= \left[\begin{array}{c}e^{i(\kappa_x^m x_i+\kappa_y^n y_i)} \\-ik_u^{(m,n)}e^{i(\kappa_x^m x_i+\kappa_y^n y_i)}\end{array}\right],
\mathbf{a} = \left[a_{mn}\right], m,n = -R, -R+1, \cdots, R-1, R.
\end{align}
\end{enumerate}
In summary, by combining Eqs. (\ref{eq:dirichlet}), (\ref{eq:qbc}), and (\ref{eq:radiation}), the whole system that is enforcing the Dirichlet boundary condition, the quasi-periodic conditions, and the upward radiation condition can be written as
\begin{align}
\left[\begin{array}{ccc}
\mathbf{A} & \mathbf{B} & \mathbf{0} \\ 
\mathbf{P} & \mathbf{Q} & \mathbf{0} \\ 
\mathbf{Z} & \mathbf{V} & \mathbf{W}
\end{array}\right]
\left[\begin{array}{c}
\mathbf{c} \\ 
\mathbf{d} \\ 
\mathbf{a}\end{array}\right] = \left[\begin{array}{c}\mathbf{f} \\\mathbf{0} \\\mathbf{0}\end{array}\right].
\end{align}

The linear system is not big for the half-space problem, but it is an overdeterminded system. A backward stable least square solver in Matlab ({\tt mldvide}) can be applied to obtain an accurate solution, or one can use a Schur complement to eliminate $\mathbf{d}$ and $\mathbf{a}$ to solve the system.

\section{Periodizing scheme for multilayered media with transmission boundary conditions}
\begin{figure}[t] %  figure placement: here, top, bottom, or page
   \centering
  \includegraphics[width=4.7in]{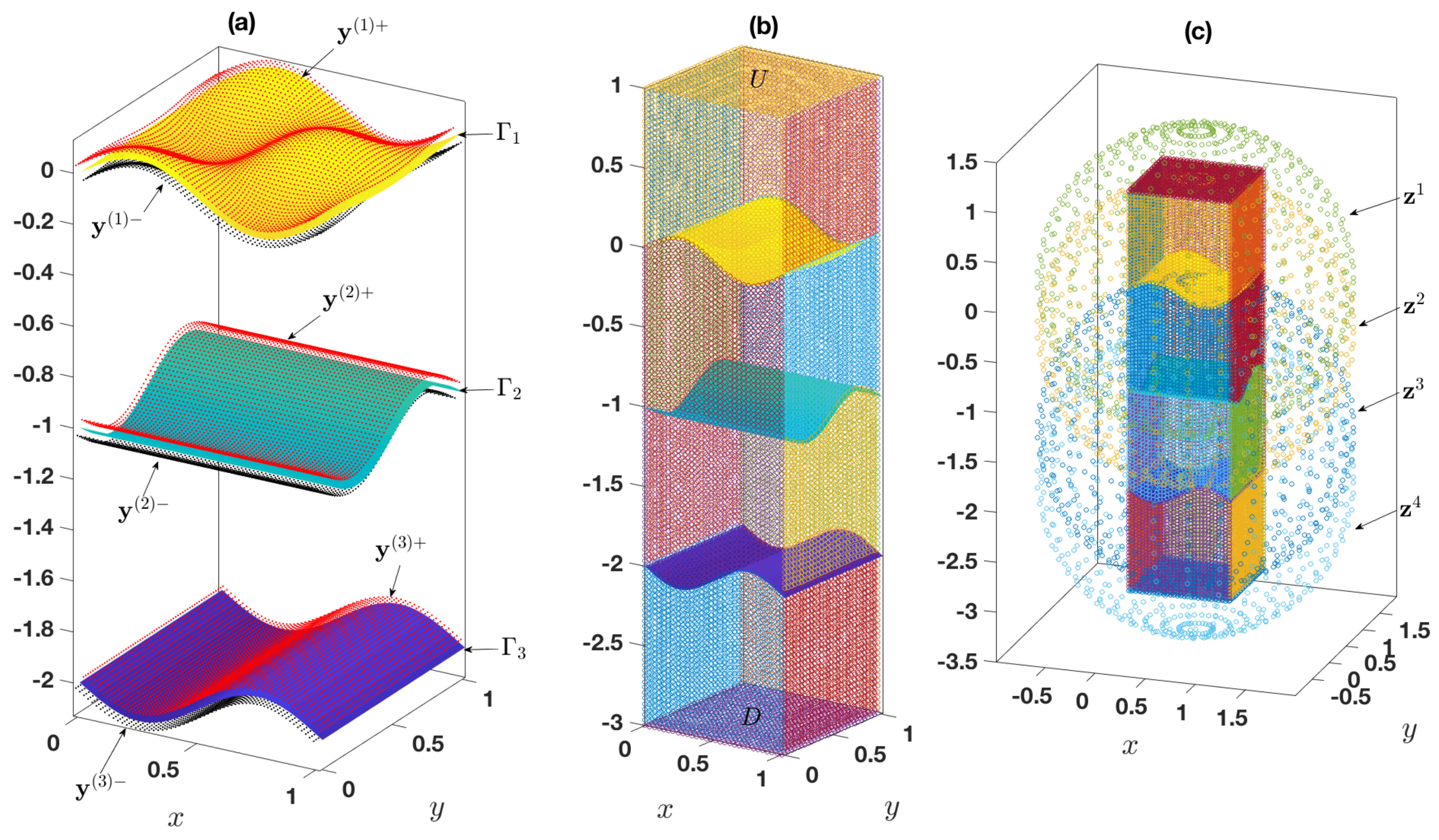}  %low resolution   
   \caption{(a) 4 layers with 3 interfaces ($\Gamma_1$, $\Gamma_2$, and $\Gamma_3$) and source points around the $i$-th interface $\mathbf{y}^{(i)-}$ (red dots) and $\mathbf{y}^{(i)+}$ (black dots) located in the normal direction below and above each interface, respectively. (b) Left, right, back, and front walls in each layer and top ($T$) and bottom ($D$) artificial layer. (c) Proxy source points (color online).}
   \label{fig:domain}
\end{figure}

In this section, multilayered media consisting of $I$ interfaces $\{\Gamma_i\}_{i=1}^I$ with transmission boundary conditions are considered. Due to the nature of the multilayered structure, it is inevitable to use many subscript and superscript indices in the notation. In each layer, the MFS representation is denoted by $\{u_i(\mathbf{x})\}_{i=1}^{I+1}$, namely,
\begin{align}
&u_1(\mathbf{x}) = \sum_{j=1}^N  \sum_{m=-1}^1 \sum_{n=-1}^1c^{(1)-}_j \alpha_x^m \alpha_y^n G_{k_1}(\mathbf{x}, \mathbf{y}^{(1)-}_j+m\mathbf{e}_x+n\mathbf{e}_y)+\sum_{j=1}^Ld^{(1)}_jG_{k_1}(\mathbf{x},\mathbf{z}^1_j), \label{eq:mfsmulti1}\\
&u_i(\mathbf{x}) =  \sum_{j=1}^N  \sum_{m=-1}^1 \sum_{n=-1}^1c^{(i)-}_j \alpha_x^m \alpha_y^n G_{k_{i}}(\mathbf{x}, \mathbf{y}^{(i)-}_j+m\mathbf{e}_x+n\mathbf{e}_y)\nonumber\\
&~~~~+ \sum_{j=1}^N  \sum_{m=-1}^1 \sum_{n=-1}^1c^{(i)+}_j \alpha_x^m \alpha_y^n G_{k_{i}}(\mathbf{x}, \mathbf{y}^{(i-1)+}_j+m\mathbf{e}_x+n\mathbf{e}_y)+\sum_{j=1}^Ld^{(i)}_jG_{k_{i}}(\mathbf{x},\mathbf{z}^{i}_j), i=2,3,\cdots I, \label{eq:mfsmulti2}\\
&u_{I+1}(\mathbf{x}) = \sum_{j=1}^N  \sum_{m=-1}^1 \sum_{n=-1}^1c^{(I+1)+}_j \alpha_x^m \alpha_y^n G_{k_{I+1}}(\mathbf{x}, \mathbf{y}^{(I)+}_j+m\mathbf{e}_x+n\mathbf{e}_y)+\sum_{j=1}^Ld^{(I+1)}_jG_{k_{I+1}}(\mathbf{x},\mathbf{z}^{I+1}_j), \label{eq:mfsmulti3}
\end{align}
where $\{\mathbf{y}^{(i)+}\}_{i=1}^{I}$ and $\{\mathbf{y}^{(i)-}\}_{i=1}^{I}$ are MFS source points placed above and below the $i$-th interface, respectively (Fig. \ref{fig:domain}(a)), and $\{\mathbf{z}^{(i)}\}_{i=1}^{I+1}$ are the proxy source points on a sphere that encloses the unit cell in the $i$-th layer (Fig. \ref{fig:domain}(c)). The incident wave is present only in the top layer. Thus, transmission boundary conditions are 
\begin{align}
\left. u_1 \right\vert_{\Gamma_1} - \left.  u_2\right\vert_{\Gamma_1} = -u^{inc} &\mbox{ and } \left. u_i\right\vert_{\Gamma_i}-\left.  u_{i+1} \right\vert_{\Gamma_i}= 0, i=2,3,\cdots I,\\
\left. \frac{\partial u_1}{\partial \mathbf{n}} \right\vert_{\Gamma_1} - \left.  \frac{\partial u_2}{\partial \mathbf{n}}\right\vert_{\Gamma_1} = -\frac{\partial u^{inc}}{\partial \mathbf{n}} &\mbox{ and } \left. \frac{\partial u_i}{\partial \mathbf{n}}\right\vert_{\Gamma_i}-\left.  \frac{\partial u_{i+1}}{\partial \mathbf{n}} \right\vert_{\Gamma_i}= 0, i=2,3,\cdots I.
\end{align}
The left, right, back and front walls in the $i$-th layer are denoted by $\Gamma_{L_i}$, $\Gamma_{R_i}$, $\Gamma_{B_i}$ and $\Gamma_{F_i}$, respectively (Fig. \ref{fig:domain}(b)). The quasi-periodic conditions must be enforced in each layer  by
\begin{align}
&\left. u \right\vert_{\Gamma_{L_i}} - \alpha_x^{-1} \left.  u \right\vert_{\Gamma_{R_i}} = 0, ~\left. \frac{\partial u}{\partial \mathbf{n}}\right\vert_{\Gamma_{L_i}} - \alpha_x^{-1} \left. \frac{\partial u}{\partial \mathbf{n}}\right\vert_{\Gamma_{R_i}} = 0\\
&\left.  u \right\vert_{\Gamma_{B_i}} - \alpha_y^{-1} \left. u\right\vert_{\Gamma_{F_i}} = 0 ,~ \left. \frac{\partial u}{\partial \mathbf{n}}\right\vert_{\Gamma_{B_i}} - \alpha_y^{-1} \left. \frac{\partial u}{\partial \mathbf{n}}\right\vert_{\Gamma_{F_i}} = 0,
\end{align}
for  $i=1,2,\cdots I+1$. Finally, unlike the Dirichlet problem, the scattered wave presents in both the top and bottom layers. Therefore, the upward and downward Rayleigh-Bloch radiation conditions
\begin{align}
u^u_{RB}(\mathbf{x}) &=\sum_{m=-\infty}^{\infty} \sum_{n = -\infty}^\infty a^u_{mn}e^{i(\kappa_x^m x+\kappa_y^n y +k_u^{(m,n)}(z-z_u)) }, z \geq z_u,\\
u^d_{RB}(\mathbf{x}) &=\sum_{m=-\infty}^{\infty} \sum_{n = -\infty}^\infty a^d_{mn}e^{i(\kappa_x^m x+\kappa_y^n y -k_d^{(m,n)}(-z+z_d)) }, z \leq z_d,
\end{align}
where $k_u^{(m,n)} = \sqrt{k_1^2-(\kappa_x^m)^2-(\kappa_y^n)^2}$ and $k_d^{(m,n)} = \sqrt{k_{I+1}^2-(\kappa_x^m)^2-(\kappa_y^n)^2}$,
have to be applied to $U = \{(x,y,z_u)| x\in [0,e_x], y\in[0, e_y]\}$ and $D = \{(x,y,z_d)|x\in [0,e_x], y\in[0, e_y] \}$, respectively. In summary, by applying transmission boundary conditions, the quasi-periodic conditions, and the upward and downward Rayleigh-Bloch radiation conditions on Eqs (\ref{eq:mfsmulti1})$\sim$(\ref{eq:mfsmulti3}), a linear system for the MFS coefficients, proxy source strengths, and the Bragg coefficients can be obtained. For simplicity's sake, a matrix structure for three interfaces ($I=3$) with four layers is presented. Let the collection of MFS coefficients, proxy coefficients, Bragg coefficients, and incident waves be
\begin{align}
\mathbf{c}& =\left[\begin{array}{cccccc}\mathbf{c}^{(1)-}  & \mathbf{c}^{(2)-} & \mathbf{c}^{(2)+} & \mathbf{c}^{(3)-} & \mathbf{c}^{(3)+} & \mathbf{c}^{(4)+}\end{array}\right]^t,\nonumber\\
\mathbf{d}&=\left[\begin{array}{cccc}\mathbf{d}^{(1)}  & \mathbf{d}^{(2)} & \mathbf{d}^{(3)} & \mathbf{d}^{(4)} \end{array}\right]^t\nonumber\\
\mathbf{a}&=\left[\begin{array}{cc}\mathbf{a}^{u}  & \mathbf{a}^{d} \end{array}\right]^t,\nonumber\\
\mathbf{f}&=\left[\begin{array}{cc}-\mathbf{u}^{inc}  & -\mathbf{\frac{\partial u}{\partial \mathbf{n}}}^{inc} \end{array}\right]^t,
\end{align}
respectively. Then, the linear system for four layers is
\begin{align}
\left[\begin{array}{cccccccccccc}
A_{1,1} & A_{1,2} & A_{1,3} & 0           & 0          &  0          & B_{1,1} & B_{1,2} & 0           & 0           & 0           & 0 \\
0           & A_{2,2} & A_{2,3} & A_{2,4} & A_{2,5} &  0          & 0          & B_{2,2}  & B_{2,3} & 0           & 0           & 0 \\
0           & 0           & 0          & A_{3,4} & A_{3,5} & A_{3,6} & 0          & 0            & B_{3,3} & B_{3,4} & 0            & 0 \\
P_{1,1} & 0           & 0          & 0           & 0          & 0           & Q_{1}   & 0            & 0           & 0           & 0           & 0 \\
0           & P_{2,2} & P_{2,3}& 0           & 0          & 0           & 0          & Q_{2}     & 0           & 0           & 0           & 0 \\
0           & 0           & 0          & P_{3,4} & P_{3,5}& 0           & 0          & 0            & Q_{3}    & 0           & 0           & 0 \\
0           & 0           & 0          & 0           & 0          & P_{4,6} & 0          & 0            & 0           & Q_{4}    & 0           & 0 \\
Z_{1}    & 0            & 0         & 0           & 0          & 0           & V_{1}    & 0           & 0            & 0           & W_{1}   & 0 \\
0           & 0           & 0          & 0           & 0          & Z_{2}    & 0           & 0           & 0            & V_{2}    & 0           & W_{2}\end{array}\right]
\left[\begin{array}{c}
\mathbf{c}^{(1)-} \\
 \mathbf{c}^{(2)-} \\
\mathbf{c}^{(2)+} \\
\mathbf{c}^{(3)-} \\
\mathbf{c}^{(3)+} \\
\mathbf{c}^{(4)+} \\
\mathbf{d}^{(1)} \\
\mathbf{d}^{(2)} \\
\mathbf{d}^{(3)} \\
\mathbf{d}^{(4)} \\
\mathbf{a}^{u}  \\
\mathbf{a}^{d} \end{array}\right] = \left[\begin{array}{c}\mathbf{f} \\0 \\0 \\0 \\0 \\0 \\0 \\0 \\0\end{array}\right]. \label{eq:matrix}
\end{align}
The derivation of all the matrix components in Eq. (\ref{eq:matrix}) is very similar to that of the Dirichlet problem in the previous section. Thus, instead of presenting a detailed formula for each element in the matrix, what each part represents is explained. $A_{i,j}$ is the interaction between target points on the layer interfaces and MFS source points. $B_{i,j}$ is the interaction between target points and proxy source points. Therefore, the first three rows represent the transmission boundary conditions on each interface. $P_{i,j}$ and $Q_i$ are the difference between the fields on the left, right, back, and front walls due to MFS source points and proxy points, respectively.  Thus, the 4$-$7th rows enforce the quasi-periodicity in each layer. $Z_1$ and $V_1$ are the interaction between the artificial layer $U$ and MFS source points and proxy points. $Z_2$ and $V_2$ are the interaction between the artificial layer $D$ and MFS source points and proxy points. Finally, $W_1$ and $W_2$ are  the Rayleigh-Bloch modes for $U$ and $D$. Thus, the last two rows apply the radiation conditions at the top and bottom layers.

Now, by rearranging unknowns as
\begin{align}
\mathbf{d}^{(1)'} &= \left[\begin{array}{c}\mathbf{d}^{(1)} \\ \mathbf{a}^u\end{array}\right],  P'_{1,1} = \left[\begin{array}{c}P_{1,1} \\Z_{1}\end{array}\right], Q'_{1} = \left[\begin{array}{cc}Q_1 & 0 \\V_1 & W_1\end{array}\right],\\
\mathbf{d}^{(4)'} &= \left[\begin{array}{c}\mathbf{d}^{(4)} \\ \mathbf{a}^d\end{array}\right], P'_{4,6} = \left[\begin{array}{c}P_{4,6} \\Z_{2}\end{array}\right], Q'_{4} = \left[\begin{array}{cc}Q_4 & 0 \\V_2 & W_2\end{array}\right],
\end{align}
the linear system simplifies to
\begin{align}
\left[\begin{array}{ccccccccccc}
A_{1,1} & A_{1,2} & A_{1,3} & 0           & 0          &  0          & B_{1,1}     & B_{1,2}   & 0               & 0                        \\
0           & A_{2,2} & A_{2,3} & A_{2,4} & A_{2,5} &  0          & 0               & B_{2,2}   & B_{2,3}      & 0                        \\
0           & 0           & 0          & A_{3,4} & A_{3,5} & A_{3,6} & 0                & 0             & B_{3,3}      & B_{3,4}              \\
P'_{1,1} & 0           & 0          & 0           & 0          & 0           & Q'_{1}        & 0             & 0                & 0                      \\
0           & P_{2,2} & P_{2,3}& 0           & 0          & 0           & 0          	      & Q_{2}       & 0               & 0                       \\
0           & 0           & 0          & P_{3,4} & P_{3,5}& 0           & 0          	      & 0              & Q_{3}        & 0                       \\
0           & 0           & 0          & 0           & 0          & P'_{4,6} & 0                & 0              & 0                & Q'_{4}                \end{array}\right]
\left[\begin{array}{c}
\mathbf{c}^{(1)-} \\
 \mathbf{c}^{(2)-} \\
\mathbf{c}^{(2)+} \\
\mathbf{c}^{(3)-} \\
\mathbf{c}^{(3)+} \\
\mathbf{c}^{(4)+} \\
\mathbf{d}^{(1)'} \\
\mathbf{d}^{(2)} \\
\mathbf{d}^{(3)} \\
\mathbf{d}^{(4)'}  \end{array}\right] = \left[\begin{array}{c}\mathbf{f} \\0 \\0 \\0  \\0 \\0 \\0\end{array}\right].
\end{align}
Then, \textcolor{brown}{the} Schur complement further reduces the system to
\begin{align}
\left[\begin{array}{ccccccccccc}
A'_{1,1} & A'_{1,2} & A'_{1,3} & 0           & 0          &  0                 \\
0           & A'_{2,2} & A'_{2,3} & A'_{2,4} & A'_{2,5} &  0            \\
0           & 0           & 0          & A'_{3,4} & A'_{3,5} & A'_{3,6}   \end{array}\right]
\left[\begin{array}{c}
\mathbf{c}^{(1)-} \\
 \mathbf{c}^{(2)-} \\
\mathbf{c}^{(2)+} \\
\mathbf{c}^{(3)-} \\
\mathbf{c}^{(3)+} \\
\mathbf{c}^{(4)+} \end{array}\right] = \left[\begin{array}{c}\mathbf{f} \\0 \\0 \end{array}\right],\label{final_system}
\end{align}
where
\begin{align}
A'_{1,1} &= A_{1,1}-B_{1,1}Q_1^{'\dag}P'_{1,1}, A'_{1,2} = A_{1,2}-B_{1,2}Q_2^{\dag}P_{2,2}, A'_{1,3} = A_{1,2}-B_{1,2}Q_2^{\dag}P_{2,3},\nonumber\\
A'_{2,2} &=A_{2,2} -B_{2,2}Q_2^{\dag}P_{2,2}  ,  A'_{2,3} =A_{2,3} -B_{2,2}Q_2^{\dag}P_{2,3} \nonumber\\
A'_{2,4} &= A_{2,4} - B_{2,3}Q_3^{\dag}P_{3,4}, A'_{2,5} = A_{2,5} - B_{2,3}Q_3^{\dag}P_{3,5} \nonumber\\
A'_{3,4} &= A_{3,4}-B_{3,3}Q_3^{\dag}P_{3,4}, A'_{3,5} = A_{3,5}-B_{3,3}Q_3^{\dag}P_{3,5}, A'_{3,6} = A_{3,6}-B_{3,4}Q_4^{'\dag}P'_{4,6}, 
\end{align}
and $Q_i^\dag$ represents the pseudo-inverse of $Q_i$. All the additional unknowns created by the periodizing method are eliminated. The MFS coefficients can be found by solving the system with the pseudo-inverse of the matrix. For multilayered media, each additional interface adds one more row in Eq. (\ref{final_system}).
%\begin{figure}[t] %  figure placement: here, top, bottom, or page
%   \centering
%   \includegraphics[width=5in]{dirichlet_field} 
%   \caption{Total field ($u+u^{inc}$) from (a) $g(x,y) = 0$ and (b) $g(x,y) = 0.2\sin{(2\pi x)}\cos{(2\pi y)}$ with $k = 15$, $\theta^{inc} = \pi/3$, $\phi^{inc} = 5\pi/6$ (color online)}
%   \label{fig:dirichlet_sol}
%\end{figure}
%\begin{figure}[t] %  figure placement: here, top, bottom, or page
%   \centering
%   \includegraphics[width=5in]{convergence_N_M_flat_Dirichlet_k_10pi} 
%   \caption{ Flux and pointwise error for flat interface with Dirichlet boundary condition, $k = 10\pi$, $\theta^{inc} = \pi/4$, $\phi^{inc} = 5\pi/6$,  $N_w = N_T = 30^2$, and $Q = -10 \cdots 10$ as a function of (a) number of source points $N = 5^2, 10^2, \cdots, 80^2$ while $M = 60\times60$ and (b) number of proxy source points $L = 5^2, 10^2, \cdots, 70^2$ while $N=50^2$ (color online).}
%   \label{fig:convergence_dirichlet}
%\end{figure}

\section{Numerical results}\label{numerical_results}
In this section, numerical examples of two-layer structures with the Dirichlet boundary condition and multilayered media with transmission conditions are presented. All computations were performed on a workstation with dual 2.6 GHz Xeon E5-2697v3 processors and 128GB RAM using Matlab R2015b. The periods in both the $x$- and $y$-axis are assumed to be 1. The MFS source points are uniformly placed at $\tau$ away from the surface in the normal direction. Target points on each surface and points on all the surrounding walls are uniformly distributed. Throughout the examples, $N$, $M$, $L$, $N_t$, $N_w$, $R$ denote the number of MFS source points, target points, proxy source points, points on top and bottom walls, points on side walls, and Rayleigh Bloch expansion terms, respectively.
\begin{figure}[h] %  figure placement: here, top, bottom, or page
   \centering
   \includegraphics[width=3in]{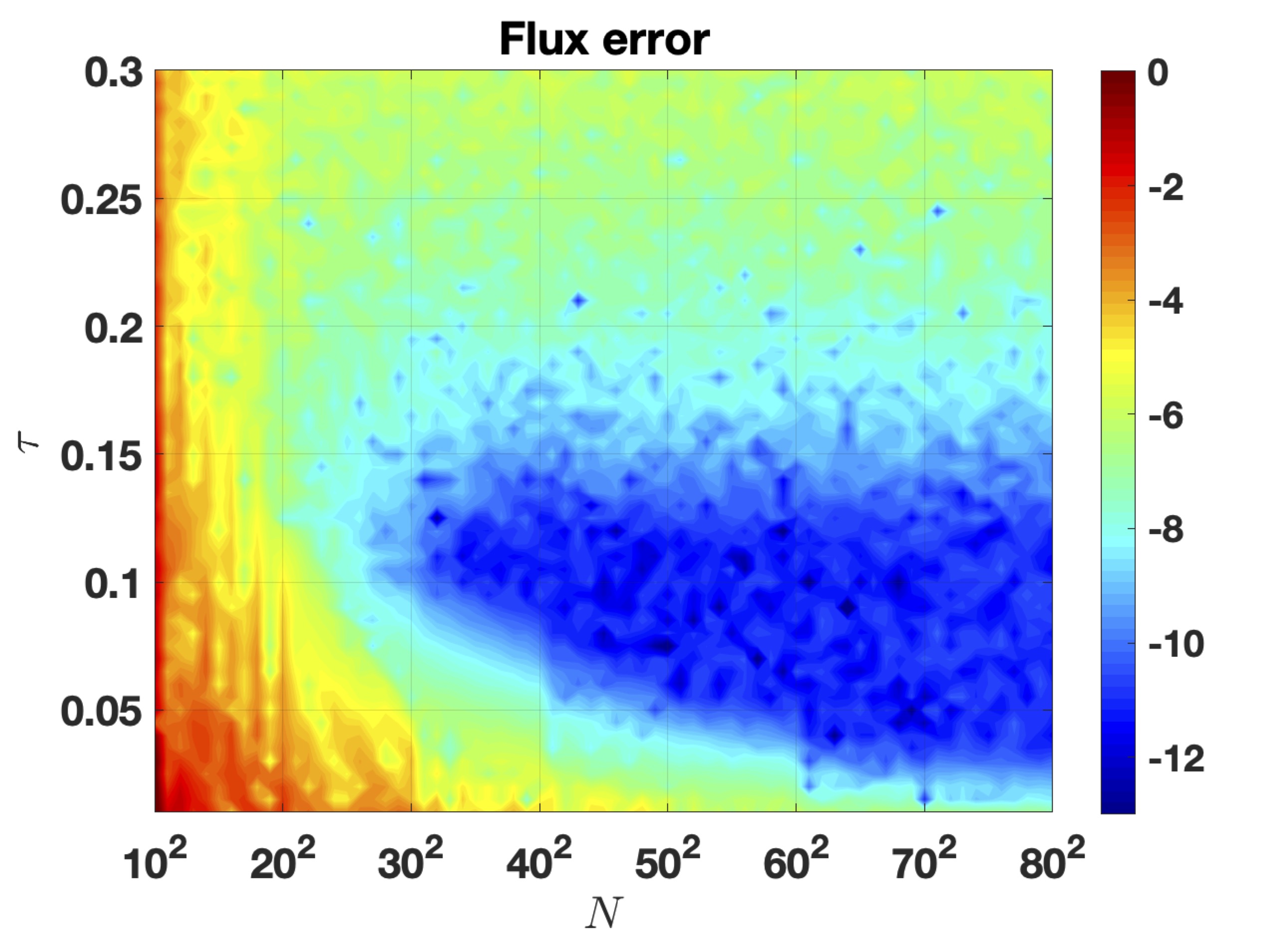} 
   \caption{$\log_{10}$(Flux error) for two-layered media with a Dirichlet boundary condition  on the interface $g(x,y) = 0.1\sin{(2\pi x)}\cos{(2\pi y)}$ with $k_1 = 10$ and incident angle $\theta^{inc} = \pi/4$ and $\phi^{inc} = 5\pi/6$ for $N=10^2, 20^2, \cdots 80^2$ and MFS source points displacement $\tau = 0.01, \cdots, 0.3$ along the normal direction (color online).}
   \label{fig:mfs_error}
\end{figure}

\begin{figure}[tt] %  figure placement: here, top, bottom, or page
   \centering
   \includegraphics[width=5.2in]{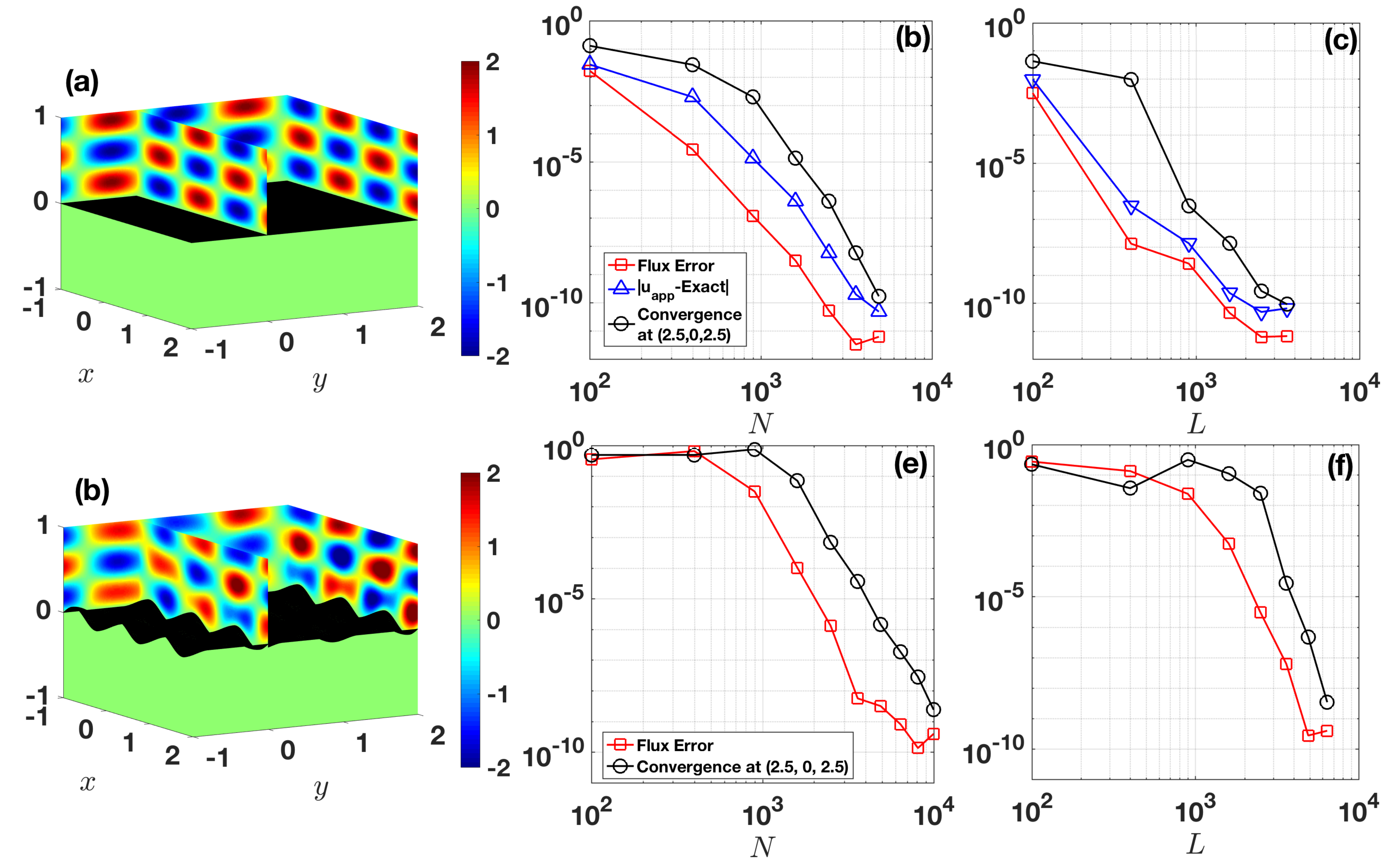} 
   \caption{Two-layered media with a Dirichlet boundary condition with $k_1 = 10$ and incident angle $\theta^{inc} = \pi/4$ and $\phi^{inc} = 5\pi/6$: (a) Total field, convergence, errors with respect to (a) $N$ and (c) $L$ for $g(x,y) = 0$. (d) Total field, convergence and flux error with respect to (e) $N$ and (f) $L$ for $g(x,y) = 0.1\sin{(2\pi x)}\cos{(2\pi y)}$   (color online).}
   \label{fig:dirichlet_data}
\end{figure}
As one measure of accuracy, the conservation of flux or energy is used, namely,
\begin{align}
\sum_{m,n} k_u^{(m,n)}|a_{mn}^u|^2+\sum_{m,n} k_d^{(m,n)}|a_{mn}^d|^2 = k_1\cos{\phi^{inc}}.
\end{align}
In other words, numerically computed energy is compared with the input energy and their relative difference is defined as a flux or an energy deficiency error:
\begin{align}
\mbox{flux error} : = \left|\frac{\sum_{m,n} k_u^{(m,n)}|a_{mn}^u|^2+\sum_{m,n} k_d^{(m,n)}|a_{mn}^d|^2-k_1\cos{\phi^{inc}}}{k_1\cos{\phi^{inc}}}\right|.
\end{align}
It has been shown that pointwise error and the flux error behave very similarly in 2D problems \cite{helmholtz_periodic}.

\begin{figure}[t] %  figure placement: here, top, bottom, or page
   \centering
   \includegraphics[width=5.2in]{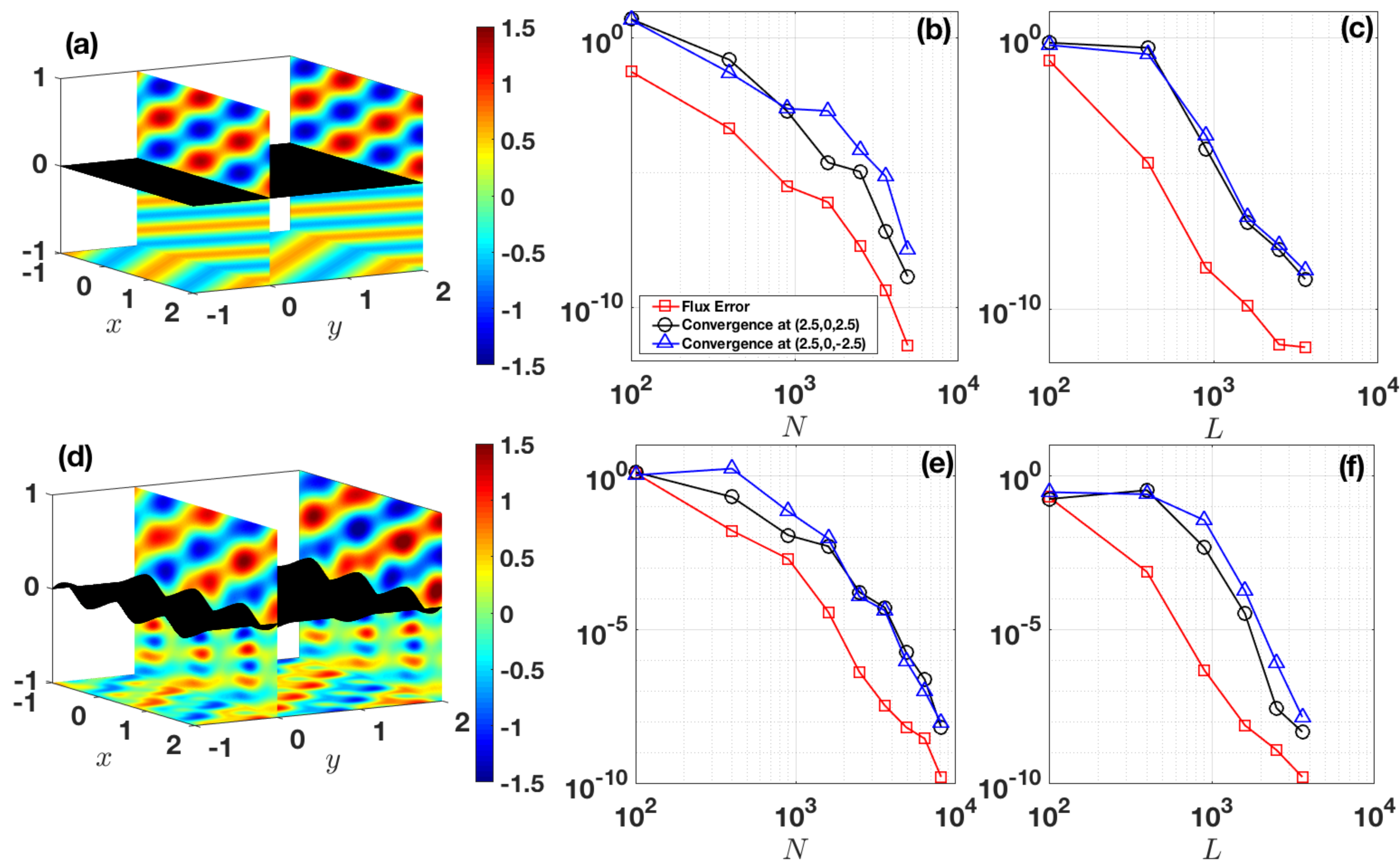} 
   \caption{Two-layered media with transmission boundary conditions with $k_1 = 10$, $k_2 = 20$ and incident angle $\theta^{inc} = \pi/4$ and $\phi^{inc} = 5\pi/6$ : (a) Total field, convergence and flux error with respect to (b) $N$ and (c) $L$ for $g(x,y) = 0$. (d) Total field, convergence and flux error with respect to (e) $N$ and (f) $L$ for $g(x,y) = 0.1\sin{(2\pi x)}\cos{(2\pi y)}$  (color online).}
   \label{fig:two_layer_data}
\end{figure}

\subsection{Two-layered media}

The flat $g(x,y) = 0$ and corrugated $g(x,y) = 0.1\sin{(2\pi x)}\cos{(2\pi y)}$ surfaces are considered and the Dirichlet boundary condition is imposed on the surface. In both examples, the wavenumber in the top layer is set to $k = 10$ and the incident angle is fixed at $\theta^{inc} = \pi/4$ and $\phi^{inc} =5\pi/6$.  In Fig. \ref{fig:mfs_error}, flux errors are computed for $N = 10^2, 20^2, \cdots 80^2$ and $\tau = 0.01 \cdots 0.3$ for the corrugated surface. A flux error of less than $10^{-10}$ can be obtained between $\tau = 0.02$ and $0.15$. Other ways of placing the MFS source points will be investigated in our future work. In the following numerical examples, $\tau = 0.03$ is used. Figures \ref{fig:dirichlet_data}(a) and \ref{fig:dirichlet_data}(d) show the total field $u+u^{inc}$ from the flat and corrugated surfaces, respectively. Fig. \ref{fig:dirichlet_data}(b) presents the flux error (red square), absolute error between the numerical solution and the exact solution at $(2.5, 0, 2.5)$ (blue triangle), and convergence at $(2.5, 0, 2.5)$ (black circle) with respect to $N$ for the flat surface. Fig. \ref{fig:dirichlet_data}(c) presents the same quantities with respect to $L$.
%For the flat surface case, flux error (red square), absolute error between the numerical solution and the exact solution at $(2.5, 0, 2.5)$ (blue triangle), and convergence at $(2.5, 0, 2.5)$ (black circle) with respect to $N$ are presented in Fig. \ref{fig:dirichlet_data}(b) and $L$ in Fig. \ref{fig:dirichlet_data}(c).
  In Fig. \ref{fig:dirichlet_data}(b),  $N$ is varied from $10^2$ to $70^2$ while all other parameters are fixed at $L = 50^2$, $N_t = 30^2$, $N_w = 30^2$. In Fig. \ref{fig:dirichlet_data}(c), $L$ is varied from $10^2$ to $60^2$ while all other parameters are fixed at $N = 70^2$, $N_t = 20^2$, $N_w = 20^2$. For the corrugated surface, $N$ is varied from $10^2$ to $100^2$ in  Fig. \ref{fig:dirichlet_data}(e) while all other parameters are fixed at $L = 80^2$, $N_t = 30^2$, $N_w = 30^2$. In Fig. \ref{fig:dirichlet_data}(c), $L$ is varied from $10^2$ to $80^2$ while all other parameters are fixed at $N = 100^2$, $N_t = 30^2$, $N_w = 30^2$. In both numerical examples, the number of target points is maintained to be $M = (1.1N)^2$ and $R = 10$ Rayleigh Bloch modes are used. Flux errors of $6.2\times10^{-12}$ and $3.9\times10^{-10}$ are obtained for the flat and corrugated surfaces, respectively. 

\begin{figure}[t] %  figure placement: here, top, bottom, or page
   \centering
      \includegraphics[width=5.3in]{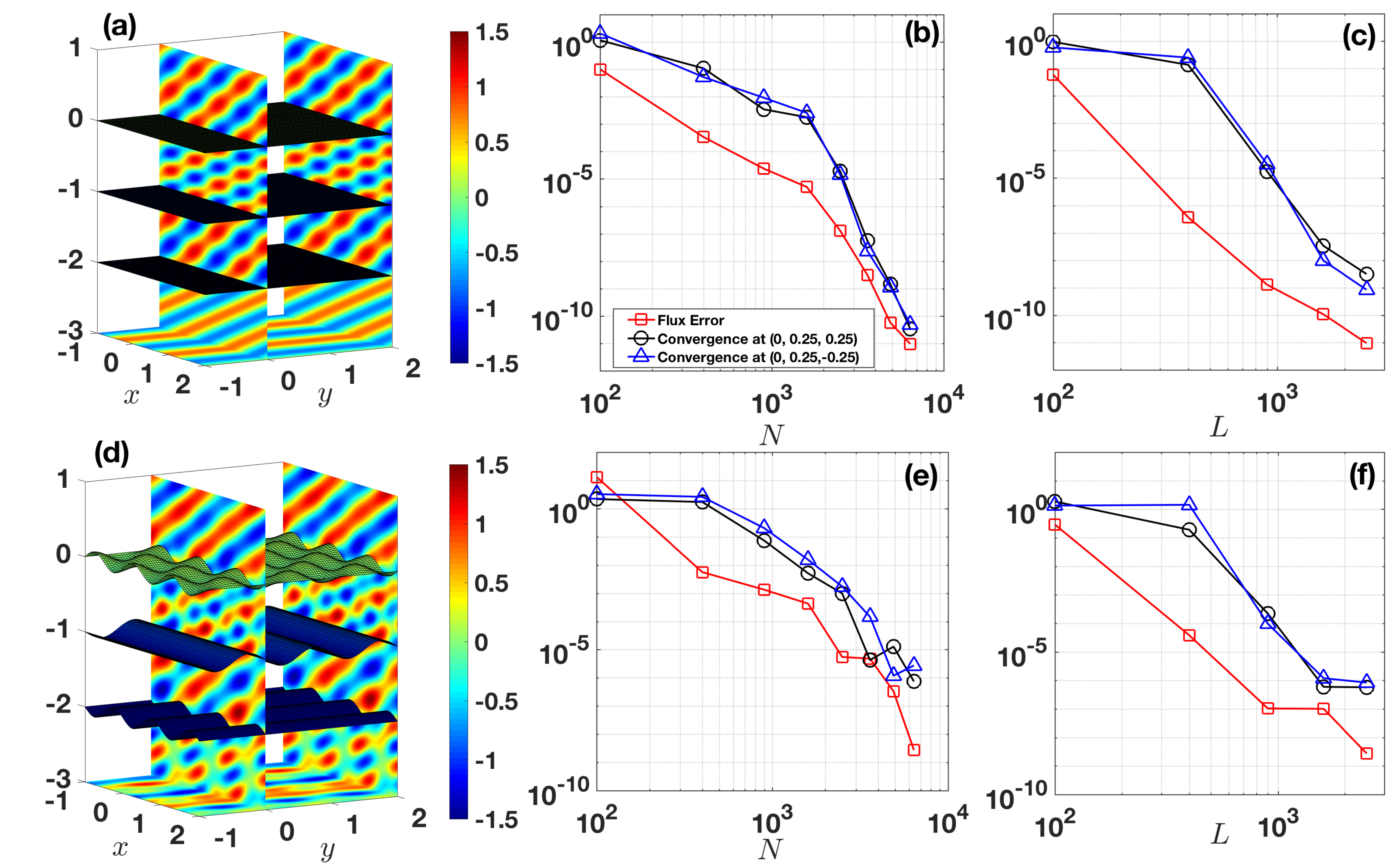} 
   \caption{Four-layered media with transmission boundary conditions with $k_1 = 3\pi$, $k_2 = 3\pi\sqrt{2}$, $k_3 = 3\pi$, $k_4 = 3\pi\sqrt{2}$ and incident angle $\theta^{inc} = 0$ and $\phi^{inc} = 5\pi/6$ : (a) Total field, convergence and flux error with respect to (b) $N$ and (c) $L$ for $g_1(x,y) = 0$, $g_2(x,y) = -1$, and $g_2(x,y) = -2$. (d) Total field, convergence and flux error with respect to (e) $N$ and (f) $L$ for $g_1(x,y) = 0.1\sin{(2\pi x)}\cos{(2\pi y)}$, $g_2(x,y) = -0.1\sin{(2\pi y)}-1$ and $g_3(x,y) = -0.1\sin{(2\pi x)}-2$ (color online).}
   \label{fig:four_layer_data}
\end{figure}
\begin{figure}[h] %  figure placement: here, top, bottom, or page
   \centering
   \includegraphics[width=5.0in]{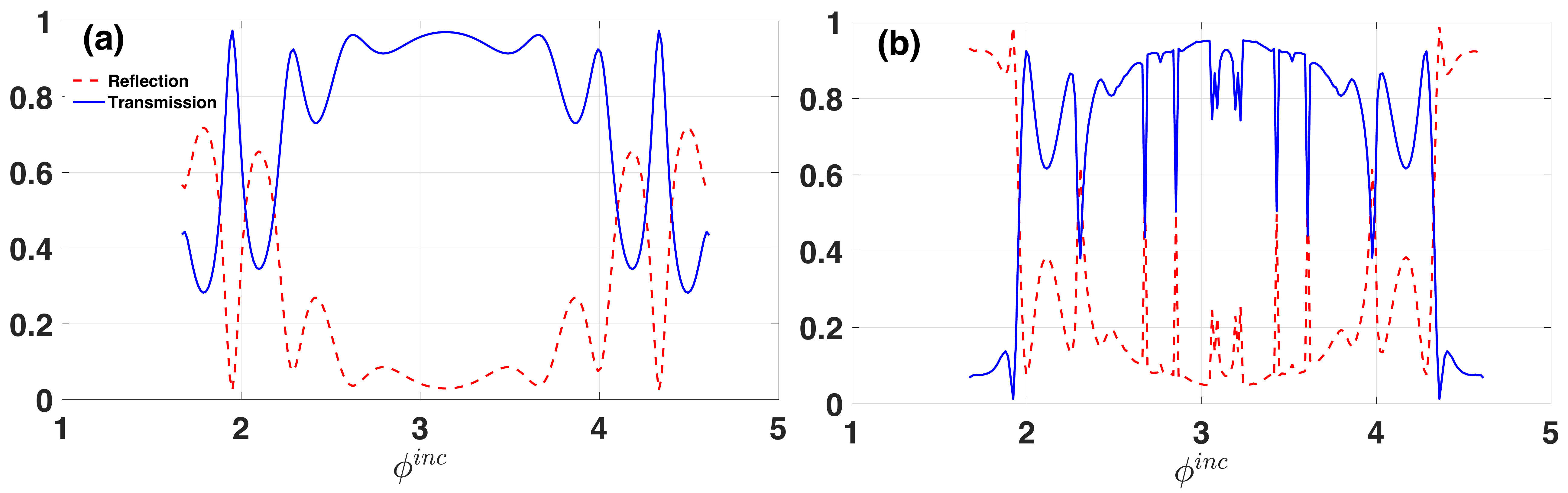} 
   \caption{Reflection (dashed red line) and transmission (solid blue line) from (a) the flat and (b) the corrugated multilayered media (color online).}
   \label{fig:reflection_data}
\end{figure}

The transmission boundary conditions are applied to the flat and corrugated surfaces in Fig. \ref{fig:two_layer_data}. In both examples, $k_1 = 10$, $k_2 = 20$, $\theta^{inc} = \pi/4$, and $\phi^{inc} = {5\pi}/{6}$ are used.  The number of points on each side wall $N_w = 30^2$, top and bottom fictitious layers $N_t = 30^2$, and  the number of Bragg mode $R = 10$ are used for all computations. The first row of Fig. \ref{fig:two_layer_data} presents numerical results for the flat surface: (a) total field, (b) convergence with respect to MFS points $N = 5^2, 10^2, \cdots,70^2$ for the fixed $L = 60^2$, and (c) convergence with respect to proxy source points $L = 10^2, 20^2,\cdots, 60^2$ while $N = 70^2$. The second row of Fig. \ref{fig:two_layer_data} presents the numerical results for the corrugated surface: (d) total field, (e) convergence with respect to MFS points $N = 5^2, 10^2, \cdots, 90^2$ for the fixed $L = 60^2$, and (f) convergence with respect to proxy source points $L = 10^2, 20^2, \cdots, 60^2$ while $N = 90^2$.  In all convergence plots, flux error (red square) and pointwise convergence at  $(2.5, 0, \pm2.5)$ (black circle and blue triangle) are displayed. Flux errors $3.9\times10^{-12}$ and $1.7\times10^{-10}$ are obtained for the flat and corrugated surfaces, respectively. Lastly, for the corrugated surface, the incident angle is chosen at the Wood anomaly in the first layer ($k_1 = 10$, $k_2 = 20$, $\theta^{inc} = \pi/4$, and $\phi^{inc} = 3.029937387$). The flux error $5.3\times 10^{-11}$ is obtained when $N = 90^2$ and  $L = 80^2$, which is no worse than that of other incident angles. Due to space limitations, the field shape is omitted.

\subsection{Multilayered media and transmission and reflection spectra}
For multilayered media, two examples consisting of 3 interfaces (4 layers) are provided in Fig. \ref{fig:four_layer_data}. In the first example, all the layer interfaces are assumed to be flat with $g_1(x,y) = 0$, $g_2(x,y) = -1$, and $g_3(x,y) = -2$ (flat multilayered medium). In the second example, layer interfaces are described by $g_1(x,y) = 0.1\sin{(2\pi x)}\cos{(2\pi y)}$, $g_2(x,y) = -0.1\sin{(2\pi y)}-1$ and $g_3(x,y) = -0.1\sin{(2\pi x)}-2$ (corrugated multilayered medium). In both examples, the wavenumber in each layer is fixed at  $k_1 = 3\pi$, $k_2 = 3\pi\sqrt{2}$, $k_3 = 3\pi$, $k_4 = 3\pi\sqrt{2}$. The incident angle is set to $\theta^{inc} = 0$ and $\phi^{inc} = {5\pi}/{6}$. The number of points on each side wall $N_w = 30^2$, top and bottom fictitious layers $N_t = 30^2$, and the number of Bragg mode $R = 10$ are used for all computations. The first row of Fig. \ref{fig:four_layer_data} presents numerical results for the flat multilayered medium: (a) total field, (b) convergence with respect to MFS points $N = 10^2, 20^2, \cdots,80^2$ for the fixed $L = 50^2$, and (c) convergence with respect to proxy source points $L = 10^2, 20^2, 20^2,\cdots, 50^2$ while $N = 80^2$. The second row of Fig. \ref{fig:two_layer_data} presents numerical results for the corrugated multilayered medium: (d) total field, (e) convergence with respect to MFS points $N = 10^2, 20^2, \cdots, 80^2$ for the fixed $L = 50^2$, and (f) convergence with respect to proxy source points $L = 10^2, 10^2, \cdots, 50^2$ while $N = 80^2$. In all convergence plots, flux error (red square) and pointwise convergences at  $(0, 2.5, 2.5)$ (black circle) and $(0, 2.5, -2.5)$ (blue triangle) are displayed. Flux errors of $9.5\times10^{-12}$ and $2.8\times10^{-9}$ are obtained for the flat and corrugated multilayered media, respectively. 

Finally, the reflection and transmission are computed for a range of incident angles for both flat and corrugated multilayered media in Fig. \ref{fig:reflection_data}.  The same geometries (four layers) and parameters are used from the previous numerical examples in Fig. \ref{fig:four_layer_data}(a) and (d).  Here, the computation is accelerated by observing that the matrix depends on incident angle only through Bloch phase or $k_x$ and $k_y$ \cite{helmholtz_periodic}. Several incident angles share the same Bloch phase. Thus, the reflection and transmission at these incident angles can be found at once. In both computations, $\theta^{inc}$ is fixed at $0$ ($k_y = 0$) and $\phi^{inc}$ is varied from $\pi/2$ to $3\pi/2$ (or $k_x$ from $-3\pi$ to $3\pi$). The reflection (red dashed line) and transmission (blue solid line) are plotted in Fig. \ref{fig:reflection_data}(a) and (b) for the flat and corrugated multilayered media, respectively.  The average flux error is maintained at $10^{-4}$. The computation is accelerated about three times (three incident angles share the same $k_x$), and taking about 12 hours to compute. Note that the computation can be further accelerated by precomputing some matrix components that are independent of the Bloch phase at the cost of computer memory.

\section{Conclusion}
A periodizing method is combined with the method of fundamental solutions for wave scattering from doubly-periodic 3D multilayered media. The scheme is robust for all scattering parameters and does not use singular quadratures at the cost of introducing artificial source points near the surfaces. Numerical examples show 9- to 10-digit accuracy at a moderate frequency region. The reflection and transmission spectra will be useful for application scientists/engineers for their studies in meta-materials, diffraction gratings, and medical imaging. Choosing an optimal location of source points is one of the drawbacks of the method. The code will be available upon request. For future work, the proposed method will be extended to include objects inside the layers, the second-kind boundary integral equation methods and/or preconditioners will be investigated for the use of an iterative or fast direct matrix solver, and the proposed method will be used for Maxwell's equations in doubly-periodic multilayered media.

\section*{Acknowledgment}

This work was supported by a grant from the Simons Foundation (\#404499, Min Hyung Cho). The author would also like to thank Dr. Alex Barnett from Flatiron Institute for helpful discussions.

%\section*{References}
\bibliography{MinHyung.bib}

\end{document}